\input amstex

\documentstyle{amsppt}

\refstyle{A}

\nologo

\hoffset .25 true in
\voffset .2 true in

\hsize=6.1 true in
\vsize=8.5 true in

\define\Adot{\bold A^\bullet}

\define\Fdot{\bold F^\bullet}
\define\Pdot{\bold P^\bullet}
\define\van{\phi_f\Fdot}
\define\dm{\operatorname{dim}}
\define\lotimes{\ {\overset L\to \otimes}\ }

\topmatter

\title A Little Microlocal Morse Theory \endtitle

\author David B. Massey \endauthor

\address{David B. Massey, Dept. of Mathematics, Northeastern University, Boston, MA, 02115, USA} \endaddress

\email{DMASSEY\@NEU.edu}\endemail


\abstract
If a complex analytic function, $f$,  has a stratified isolated critical point, then it is known that the cohomology of the Milnor fibre
of $f$ has a direct sum decomposition in terms of the normal Morse data to the strata. We use microlocal Morse theory to obtain the same
result under the weakened hypothesis that the vanishing cycles along $f$ have isolated support. We also investigate an index-theoretic
proof of this fact. 
\endabstract
\endtopmatter

\document

\noindent\S0. {\bf Introduction}  

\vskip .2in

Our goal in this paper is very modest. Let $X$ be a complex analytic space, let $\Cal S:=\{S_\alpha\}$ be a Whitney stratification of $X$
with connected strata, let
$f:X\rightarrow\Bbb C$ be a complex analytic function, let
$R$ be a principal ideal domain (p.i.d.), and let
$\Fdot$ be a bounded complex of sheaves of $R$-modules on
$X$, which is constructible with respect to $\Cal S$. (Typically, one choses the base ring $R$ to be $\Bbb Z$, $\Bbb Q$, or $\Bbb C$.) We
wish to show: if
$\bold x$ is an isolated point in the support of the vanishing cycles
$\van$, then the stalk cohomology of
$\van$ at $\bold x$ is a direct sum, with multiplicities, of the $\Fdot$-Morse modules of strata of
$X$  (i.e., hypercohomologies of normal data of strata, with coefficients in
$\Fdot$).

As we are interested in activity at an isolated point, our question is entirely local. Hence, there is no loss of generality in assuming
that $X$ is embedded in an open subset, $\Cal U$, of some affine space, that $\bold x = \bold 0$, and that $\hat f$ is an analytic
extension of
$f$ to all of
$\Cal U$.

To state our result precisely, we need to introduce terminology that will appear in Definition 3.3: a stratum $S_\alpha$ is {\it
$\Fdot$-visible} if its $\Fdot$-Morse modules are not all trivial. In Section 5, we prove:

\vskip .1in

\noindent{\bf Main Theorem} (5.3). {\it  Suppose
that $\bold 0$ is an isolated point in the support of $\phi_f\Fdot$. 

Then,  for every $\Fdot$-visible stratum $S_\alpha$, $(\bold 0, d_\bold 0\hat f)$ is an isolated point in
$\overline{T^*_{{}_{S_\alpha}}\Cal U}\cap\operatorname{im} d\hat f$, and if we let 
$k_\alpha$ equal the intersection multiplicity $\left(\overline{T^*_{{}_{S_\alpha}}\Cal U}\ \cdot\ \operatorname{im} d\hat
f\right)_{(\bold 0, d_\bold 0\hat f)}$,
then there is a (non-canonical) isomorphism
 $$ H^{i-1}(\phi_f\Fdot)_\bold 0\ \cong \ \bigoplus\Sb\Fdot{\text-visible}\\ S_\alpha\endSb\big(R^{{k_\alpha}}\otimes_{{}_R}\Bbb
H^{i-d_\alpha}(\Bbb N_\alpha,
\Bbb L_\alpha \ ; \ \Fdot)\big).$$ 
}

\vskip .1in

The knowledgeable reader may believe that this result is already known; on the cohomological level, it is a generalization of the
fundamental theorem of stratified Morse Theory: that local Morse data is the product of tangential and normal Morse data. Moreover, there
are similar results of Siersma [{\bf Si}], Tib\u ar [{\bf Ti}], and in our own work in [{\bf M2}]. However, in each of these papers, the
hypothesis is that
$f$ has a {\bf stratified} isolated critical point $\bold x$. While this assumption certainly implies that $\bold x$ is an
isolated point in
$\operatorname{supp}\van$, it is easily seen to be a strictly stronger condition.  From a philosophical point of
view, one would like to be able to conclude something about the structure of $\van$ from a condition on
$\operatorname{supp}\van$.

\vskip .2in

We first prove, in Section 2, a weak version of the main theorem; we accomplish this by using a vanishing index theorem, proved independently by
L\^e [{\bf L1}], Sabbah [{\bf Sa}], and Ginsburg [{\bf G}], combined with the tool of perverse cohomology (see [{\bf BBD}] and [{\bf K-S3}]).
Perverse cohomology allows us to take this index result about the Euler characteristic and use it to extract the individual Betti numbers. This
proof is very short. However, it has the disadvantage of being too magical; in no way does one ``see'' the
contributions from the critical points of a small perturbation of $f$. 

As a consequence of its mystical nature, the proof does not generalize to other situations where one might hope to
improve Morse-theoretic results. Moreover, one can not quite recover the integral cohomology by
this method. Hence, we wish to go back to the Morse theory proof and try to make it work.

\vskip .2in

How difficult is it to pass from the Morse theory proof in the stratified isolated case to the case of isolated
support for the vanishing cycles? Fairly, and the difficulties are subtle. Essentially, the problem is to show that $\Fdot$-invisible
strata are not relevant. Proving that such invisible strata do not matter either requires some very tedious arguments about types of
neighborhoods that one may use to calculate the normal data at points on the boundary of a small closed ball, or requires one to use the
microlocal theory developed by Kashiwara and Schapira in [{\bf K-S1}], [{\bf K-S2}], and [{\bf K-S3}].

Kashiwara and Schapira indicate ([{\bf K-S3}],
Chapt. V Notes) that the micro-support can be viewed as a generalization of stratified Morse theory -- for it applies even
to complexes of sheaves which are not constructible. In this paper, we wish to expound on two further advantages of the
microlocal approach. The first is that, by definition, the micro-support ignores superfluous strata. The second is that
the general machinery of the microlocal theory enables one to circumvent many Morse theoretic arguments involving choices
of special kinds of neighborhoods. 

However, we still find the language and ideas of stratified Morse theory to be very intuitive. Thus, our general approach
is a mixture of stratified Morse theory and the microlocal theory -- we think in terms of Morse theory, but write down
the microlocal proofs.

\vskip .4in

\noindent A brief outline of this paper is as follows. 

\vskip .1in

In Section 1, we establish some general notation and terminology that we shall use throughout the rest of the paper. Despite the fact
that the main theorem (5.3) of this paper is a statement in the complex analytic setting, we will need to work in the real subanalytic
category; hence, our notation must be more general than that which we have used in the introduction. Also in Section 1, we state and
sketch the proof of our prior result from [{\bf M2}]; we do this to provide background and because the proof of the main theorem of
this paper follows along the same lines, once we have explained the requisite microlocal theory.

In Section 2, we provide an index-theoretic proof of the main theorem in the case where the base ring is a field. In other words, we show
that there is an equality of Betti numbers. This proof is short, but not very geometrically enlightening. Moreover, even switching the
base ring to finite fields does not enable one to recover the main theorem with integral coefficients.

In Section 3, we begin our discussion of the technical details necessary to turn the sketch of a proof in Section 1 into a rigorous proof
of Theorem 5.3. We discuss switching base rings, $\Fdot$-visible strata, and the support of the vanishing cycles. The fundamental results
of this section are Proposition 3.2 (Continuity of Vanishing Support) and Theorem 3.4, which provides a conormal characterization of the
support of the vanishing cycles.

In Section 4, we summarize the basics of the microlocal theory that we shall need. For the proofs of most of the results of this section,
the reader is referred to [{\bf K-S3}]. However, for lack of any convenient reference, we provide a proofs of Theorem 4.13 and Corollary
4.15; these results relate the micro-support of a complex of sheaves to visible strata and the vanishing cycles.

Finally, in Section 5, we use the microlocal theory of Section 4 to turn our earlier proof, sketched in Section 1, into a proof of the
main theorem.

\vskip .4in

\noindent\S1. {\bf Notation and a Known Result}  

\vskip .2in

In this section, we will fix some notation from stratified Morse theory  with coefficients in
complexes of sheaves (see [{\bf G-M2}], 6.A), and then 
describe our result from [{\bf M2}], Theorem 3.2, which is based on the argument in [{\bf Si}].  We also provide a
quick sketch of the proof; the point of describing this proof is that the microlocal argument proceeds
along the same lines.

\vskip .2in

Let $R$ be a regular Noetherian ring with finite Krull dimension (e.g., $\Bbb Z, \Bbb Z/p,\ 
\Bbb Q, \text{or}\ \Bbb C$). Let $X$ be a real subanalytic set, and let $D^b(X)$ denote the derived category of
bounded complexes of
$R$-modules on
$X$. Let $\Fdot\in D^b(X)$. If there exists a real subanalytic Whitney stratification, $\Cal S$,
of $X$ with respect to which
$\Fdot$ is constructible, then we write $\Fdot\in D^b_{{}_\Cal S}(X)$ and say that $\Fdot$ is $\Bbb
R$-constructible; note that our assumptions about
$R$ guarantee that such an $\Fdot$ is perfect (see  [{\bf G-M1}], 1.4 and [{\bf K-S}], 8.4.3).   We denote the full
subcategory of
$D^b(X)$ of $\Bbb R$-constructible complexes by  $D^b_{{}_\Bbb R}(X)$.

If $X$ is complex analytic, $\Cal S$ is a complex analytic Whitney stratification of $X$, and $\Fdot\in D^b_{{}_\Cal
S}(X)$, then we say that $\Fdot$ is $\Bbb C$-constructible. We denote the full
subcategory of
$D^b(X)$ of $\Bbb C$-constructible complexes by  $D^b_{{}_\Bbb C}(X)$.

\vskip .2in

Now suppose that $X$ is a complex analytic set, $\Cal S$ is a
 complex analytic Whitney stratification of $X$, and $\Fdot\in D^b_{{}_\Cal S}(X)$. For each stratum $S_\alpha\in\Cal
S$, there exists {\it normal data to $S_\alpha$}, which is  a  pair $(\Bbb N_\alpha, \Bbb L_\alpha)$ consisting of a normal slice and
the complex link of $S_\alpha$  to $S_\alpha$; the hypercohomology modules $\Bbb H^*(\Bbb N_\alpha, \Bbb L_\alpha; \
\Fdot)$ are well-defined, and we refer to them as the {\it
$\Fdot$-Morse modules of $S_\alpha$}.

\vskip .4in

Theorem 3.2 of [{\bf M2}] is

\vskip .3in

\noindent{\bf Theorem 1.1}. {\it  Let $X$ be a complex analytic space
embedded in an open subset
$\Cal U$ of complex affine space, let $\Cal S:=\{S_\alpha\}$ be a complex Whitney stratification of $X$ with connected strata, let
$\Fdot\in D^b_{{}_\Cal S}(X)$, and let $\hat f:\Cal U\rightarrow\Bbb C$ be complex analytic. Let $f:= \hat f_{|_{X}}$, and suppose
that $f: (X, \bold 0)\rightarrow(\Bbb C, 0)$ has a stratified isolated critical point at $\bold 0$.

Then, there exists a unique set, $\{k_\alpha\}$, of non-negative integers such that,
for all
 bounded complexes, $\Fdot$, of $R$-modules on $X$ which are constructible with
respect to $\Cal S$, and, for all $i$,  $$ H^{i-1}(\phi_f\Fdot)_\bold 0 \ \cong \
\Bbb H^i(B_\epsilon \cap X, F_{f, \bold 0} \ ; \Fdot ) \ \cong \ \Bbb H^i(B_\epsilon
\cap X\cap f^{-1}(\overset\circ\to{\Bbb  D}_\eta),  B_\epsilon \cap X\cap f^{-1}(v) \
; \Fdot ) $$ $$\cong \  \bigoplus_\alpha
\big(R^{{k_\alpha}}\otimes_{{}_R}\Bbb
H^{i-d_\alpha}(\Bbb N_\alpha,
\Bbb L_\alpha \ ; \ \Fdot)\big).,$$ where:

\vskip .1in

\noindent $d_\alpha = \dm  S_\alpha$;

\vskip .1in

\noindent $B_\epsilon$ is a sufficiently small closed ball of radius $\epsilon$
centered at the origin in the ambient affine space;

\vskip .1in

\noindent $\overset\circ\to{\Bbb D}_\eta$ is an open disc of radius $\eta$, 
$0<\eta\ll \epsilon$, centered at the origin;

\vskip .1in

\noindent $v\in \overset\circ\to{\Bbb D}_\eta^*$;

\vskip .1in

\noindent $F_{f, \bold 0}$ is the Milnor fibre of $f$ at $\bold 0$ or, more 
precisely,  $F_{f, \bold 0} = B_\epsilon\cap X\cap f^{-1}(v)$.

\vskip .1in

Moreover, $k_\alpha = 1$ if $S_\alpha = \{\bold 0\}$ and, if $S_\alpha\neq\{\bold
0\}$, then for a generic choice of complex linear forms,  $L$, on the ambient affine space, for all strata $S_\alpha$, the relative
polar curve
$\Gamma_{{}_{f,L}}(S_\alpha)$ is one-dimensional (or  empty) at the origin
and  $$k_\alpha = \big(\Gamma^1_{{}_{f, L}}(S_\alpha)\cdot V(f)\big)_\bold
0  -  \big(\Gamma^1_{{}_{f, L}}(S_\alpha)\cdot V(L)\big)_\bold 0 \ =\ \left(\overline{T^*_{{}_{S_\alpha}}\Cal
U}\ \cdot\ \operatorname{im} d\hat f\right)_{(\bold 0, d_{\bold 0}\hat f)}.$$  The integer $k_\alpha$ is precisely the number of
non-degenerate critical points of a small perturbation of $f$ by $L$ which occur near the origin on the stratum
$S_\alpha$ ; more precisely, for all sufficiently small $\delta >
0$, for all complex $t$ such that $0< |t|\ll \delta$, $k_\alpha$ equals the number
of critical points of  $f+tL$  in
$\overset\circ\to{B_\epsilon}\cap S_\alpha$.
}

\vskip .3in

\noindent{\bf Sketch of the proof}:

\vskip .2in

Fix a small ball $B_\epsilon$. Since $f$ has a stratified isolated critical point at the origin, a small perturbation
of $f$, $g:=f+tL$, will have a finite number of isolated Morse stratified critical points in the interior of
$B_\epsilon$ and no critical points on the boundary. Morover, the Milnor fibre $B_\epsilon\cap V(f-v)$ will be
stratified homeomorphic to $B_\epsilon\cap V(g-v)$. Now, apply stratified Morse theory to the function
$h:=||g-v||^2$; $h^{-1}(0)$ begins at the Milnor fibre and grows out to having hypercohomology equal to the stalk
cohomology of $\Fdot$ at the origin.  The cohomology of the level sets of $h$ ``jumps'' by precisely
$\Bbb H^*(\Bbb N_\alpha,\Bbb L_\alpha; \ \Fdot)$ at each critical point on $S_\alpha$. The proof would be finished
except that one has to prove that there is no cancellation among the contributions from the various critical points.

To show that we get the direct sum decomposition in the theorem, one considers the stratified critical {\bf values}
of $g$, which can be assumed to be distinct. Thus, one can view the whole situation ``downstairs'' in a small
complex disk, $\overset\circ\to{\Bbb D}_\eta$, around the origin; one looks at the inverse image of this disk by
$g$, modulo the inverse image of
$v$, and considers what happens at the finite set of critical values. By homotoping and excising, one shows that
$\Bbb H^*\big(g^{-1}\big(\overset\circ\to{\Bbb D}_\eta\big), g^{-1}(v);\ \Fdot\big)$ breaks up as a direct sum of
the hypercohomologies of the inverse images under $g$ of a collection of small disks around the critical values
modulo other points inside the small disks. The desired conclusion follows.\qed

\vskip .4in

\noindent\S2. {\bf A Weak Version of the Main Theorem via an Index Theorem}  

\vskip .2in

In Theorem 2.4 of this section, we show how the Euler characteristic information given by the vanishing index theorem of  L\^e [{\bf L1}],
Sabbah [{\bf Sa}], and Ginsburg [{\bf G}] (see Theorem 2.1) can actually be used to find the individual Betti numbers; all that is required
is a quick application of perverse cohomology. Thus, in this section, we give a short, elegant proof of the main theorem, {\bf if one is
willing to ignore torsion}. The other sections of this paper are present precisely because, in general, we consider it unsatisfactory to
ignore torsion, and because the index-theoretic approach is not easily adaptable to other questions.

\vskip .1in

Throughout this section, we assume that our base ring, $R$, is a principal ideal domain, that $X$ is a complex analytic space
embedded in an open subset
$\Cal U$ of complex affine space, that $\Cal S:=\{S_\alpha\}$ is a complex Whitney stratification of $X$ with connected strata, that
$\Fdot\in D^b_{{}_\Cal S}(X)$, and that $\hat f:\Cal U\rightarrow\Bbb C$ is complex analytic. Recall that $d_\alpha:=\dm S_\alpha$.
Let
$f:=
\hat f_{|_{X}}$.

\vskip .2in

As $R$ is a p.i.d., the rank of a finitely-generated $R$-module is well-defined, and the Euler characteristic -- the
alternating sum of the ranks -- is a well-defined, additive function on complexes of finitely-generated $R$-modules. Recall that the
{\it characteristic cycle, 
$\operatorname{Ch}(\Fdot)$, of
$\Fdot$} in
$T^*\Cal U$ is the linear combination 
$\sum_{\alpha} m_{\alpha}(\Fdot) \left[\overline{T^*_{{}_{S_\alpha}}\Cal U}\right]$, where   the $m_\alpha(\Fdot)$
are integers determined by the Euler characteristic, as follows
$$m_\alpha(\Fdot)\ :=  \ (-1)^{\dm  X-d_\alpha}\chi\big(\Bbb H^*(\Bbb N_\alpha, \Bbb L_\alpha;\Fdot)\big)\ =\ 
(-1)^{\dm  X-d_\alpha-1}\chi(\phi_{L_{|_{\Bbb N_\alpha}}}{\bold F^\bullet}_{|_{\Bbb N_\alpha}})_{\bold x}$$  for
any point $\bold x$ in 
$S_\alpha$, with normal slice $\Bbb N_\alpha$  at $\bold x$,  and any  $L: (\Cal U, x) \rightarrow\ (\Bbb C,0)$
such that  $d_{\bold x}L$  is a non-degenerate covector at $\bold x$ (with respect to our fixed stratification; see
[{\bf G-M2}])  and
$L_{|_{S_\alpha}}$ has a Morse singularity at $\bold x$.  This cycle is independent of all the choices made (see, for
instance, [{\bf K-S3}], Chapter IX).

\vskip .3in

Let $\operatorname{Ch}(\Fdot)$ denote the characteristic cycle of $\Fdot$ in $T^*\Cal U$. Then, the result of  L\^e [{\bf L1}], Sabbah
[{\bf S}], and Ginsburg [{\bf G}] is:

\vskip .4in

\noindent{\bf Theorem 2.1}. {\it Suppose that $(\bold p, d_{\bold p}\hat f)$ is an isolated point in the intersection
$|\operatorname{Ch}(\Fdot)|\cap \operatorname{im}d\hat f$, and that $f(\bold p)=0$. Then, the Euler characteristic of the stalk
cohomology of the vanishing cycles of
$f$ is related to the intersection multiplicity of $\operatorname{Ch}(\Fdot)$ and image of $d\hat f$ by
$$\chi(\phi_{f}\Fdot)_\bold p= (-1)^{{}^{\dm X-1}}\big(\operatorname{Ch}(\Fdot)\cdot
\operatorname{im}d\hat f\big)_{(\bold p, d_{\bold p}\hat f)}.$$ }

\vskip .4in

\noindent{\it Remark 2.2}.  The proofs of L\^e, Sabbah, and Ginsburg all use the language or techniques of $\Cal D$-modules; hence,
they all use the complex field for coefficients. L\^e's Morse-theoretic proof certainly goes through without change if the base
ring is a p.i.d. Moreover, our generalization of 2.1, which appears in [{\bf M1}], Theorem 2.10 also has a purely
Morse-theoretic proof and, hence, is clearly true even when the base ring is a p.i.d.

\vskip .3in

We wish to show how one can use perverse cohomology to extract Betti number information from Theorem 2.1. We  list some properties
of the perverse cohomology and of vanishing cycles that we will need later.  The reader is referred to [{\bf BBD}] and [{\bf K-S3}].

\bigskip

The perverse cohomology functor (using middle perversity, $\mu$) on $X$, 
${}^{\mu}\negmedspace H^0$, is a functor from $D^b_{{}_\Bbb C}(X)$  to the Abelian
category of perverse sheaves on $X$. One lets ${}^{\mu}\negmedspace H^i(\Fdot)$ denote ${}^{\mu}\negmedspace H^0(\Fdot[i])$.

If $\Fdot$ is constructible with respect to $\Cal S$, then ${}^{\mu}\negmedspace H^0(\Fdot)$ is also constructible with respect to
$\Cal S$, and $\big({}^{\mu}\negmedspace H^0(\Fdot)\big)_{|_{\Bbb N_\alpha}}[-d_\alpha]$ is naturally isomorphic to
${}^{\mu}\negmedspace H^0(\Fdot_{|_{\Bbb N_\alpha}}[-d_\alpha])$.

The functor 
${}^{\mu}\negmedspace H^0$, applied to a perverse sheaf $\Pdot$ is canonically isomorphic to $\Pdot$. In addition, a bounded,
constructible complex of sheaves $\Fdot$ is perverse if and only ${}^{\mu}\negmedspace H^k(\Fdot)=0$ for all $k\neq 0$. In
particular, if $X$ is a local complete intersection, then ${}^{\mu}\negmedspace H^{{\dm X}}(\Bbb Z^\bullet_X)\cong \Bbb
Z^\bullet_X[\dm X]$ and
${}^{\mu}\negmedspace H^k(\Bbb Z^\bullet_X) = 0$ if $k\neq \dm X$.

The functor 
${}^{\mu}\negmedspace H^0$ commutes with vanishing cycles with a shift of $-1$, nearby cycles with a shift of $-1$, and 
Verdier dualizing.  That is, there are natural isomorphisms
$${}^{\mu}\negmedspace H^0 \circ \phi_f[-1] \cong
\phi_f[-1] \circ {}^{\mu}\negmedspace H^0,\hskip .2in {}^{\mu}\negmedspace H^0 \circ \psi_f[-1] \cong
\psi_f[-1] \circ {}^{\mu}\negmedspace H^0, \hskip .1in\text{and }\Cal D\circ {}^{\mu}\negmedspace H^0\cong {}^{\mu}\negmedspace H^0\circ\Cal D.$$ 

Let $\bold F^\bullet$  be a bounded  complex of sheaves on $X$ which is constructible with respect to $\Cal S$. Let
$S_{\operatorname{max}}$ be a maximal stratum (i.e., one not contained in the closure of another) which is contained in the support
of $\Fdot$, and let 
$m=\dm S_{\operatorname{max}}$. Then, $\left({}^{\mu}\negmedspace H^0(\Fdot)\right)_{|_{S_{\operatorname{max}}}}$ is isomorphic (in
the derived category) to the complex which has $\left(\bold H^{-m}(\Fdot)\right)_{|_{S_{\operatorname{max}}}}$ in degree
$-m$ and zero in all other degrees. 

In particular, $\operatorname{supp}\Fdot = \bigcup_i \operatorname{supp}{}^{\mu}\negmedspace H^i(\Fdot)$, and if $\Fdot$
is supported on an isolated point, $\bold q$, then
$H^0({}^{\mu}\negmedspace H^0 (\bold F^\bullet))_{\bold q}
\cong H^0(\bold F^\bullet)_{\bold q}.$

\vskip .4in

\noindent{\bf Lemma 2.3}. {The characteristic cycle of the perverse cohomology of $\Fdot[i]$ is given by
$$\operatorname{Ch}({}^{\mu}\negmedspace H^i(\Fdot))=(-1)^{\dm  X}\sum_\alpha b_{{}_{i-d_\alpha}}(\Bbb N_\alpha, \Bbb
L_\alpha;\Fdot)\left[\overline{T^*_{{}_{S_\alpha}}\Cal U}\right],$$ 
where $b_{{}_{j}}(\Bbb N_\alpha, \Bbb
L_\alpha;\Fdot)$ denotes the $j$-th relative Betti number, i.e., $b_{{}_{j}}(\Bbb N_\alpha, \Bbb
L_\alpha;\Fdot)={\operatorname{rk}}_{{}_R}\Bbb H^j(\Bbb N_\alpha, \Bbb L_\alpha;\Fdot)$.

}

\vskip .3in

\noindent{\it Proof}. 
$$m_\alpha({}^{\mu}\negmedspace
H^0(\Fdot[i])) =  (-1)^{\dm  X-d_\alpha-1}\chi(\phi_{L_{|_{\Bbb
N_\alpha}}}{{}^{\mu}\negmedspace
H^0(\Fdot[i])}_{|_{\Bbb N_\alpha}})_{\bold x}=$$
$$(-1)^{\dm  X}\chi\big(\phi_{L_{|_{\Bbb N_\alpha}}}[-1]{}^{\mu}\negmedspace H^0(\Fdot[i])_{|_{\Bbb
N_\alpha}}[-d_\alpha]\big)_{\bold x} = (-1)^{\dm X}\chi\big(\phi_{L_{|_{\Bbb N_\alpha}}}[-1]{}^{\mu}\negmedspace
H^0(\Fdot_{|_{\Bbb N_\alpha}}[i-d_\alpha])_{\bold x}=$$
$$
 (-1)^{\dm  X}\operatorname{rk}H^0\big(\phi_{L_{|_{\Bbb
N_\alpha}}}[-1](\Fdot_{|_{\Bbb N_\alpha}}[i-d_\alpha])\big)_{\bold x}=(-1)^{\dm  X}b_{{}_{i-d_\alpha}}(\Bbb N_\alpha, \Bbb
L_\alpha;\Fdot).\qed
$$

\vskip .4in

We can now prove the fundamental result of this section.

\vskip .4in

\noindent{\bf Theorem 2.4}. {\it Suppose that $(\bold p, d_{\bold p}\hat f)$ is an isolated point in the intersection
$|\operatorname{Ch}(\Fdot)|\cap \operatorname{im}d\hat f$, and that $f(\bold p)=0$. Then, the Betti numbers of the stalk cohomology
of the vanishing cycles of
$f$ are given by
$$b_{i-1}(\phi_{f}\Fdot)_\bold p= \sum_\alpha b_{{}_{i-d_\alpha}}(\Bbb N_\alpha, \Bbb
L_\alpha;\Fdot)\big(\overline{T^*_{{}_{S_\alpha}}\Cal U}\ \cdot\ \operatorname{im}d\hat f\big)_{(\bold p, d_\bold p\hat f)}=$$ 
$$\sum_\alpha b_{{}_{i-d_\alpha}}(\Bbb N_\alpha, \Bbb
L_\alpha;\Fdot)\left[\big(\Gamma^1_{{}_{f, L}}(S_\alpha)\cdot V(f)\big)_\bold
p  -  \big(\Gamma^1_{{}_{f, L}}(S_\alpha)\cdot V(L)\big)_\bold p\right].$$ 
}

\vskip .3in

\noindent{\it Proof}. We use the properties of perverse cohomology:

$$
b_{i-1}(\phi_{f}\Fdot)_\bold p= \operatorname{rk} H^0(\phi_f\Fdot[i-1])_\bold p = \operatorname{rk} H^0\big({}^{\mu}\negmedspace
H^0(\phi_f\Fdot[i-1])\big)_\bold p = 
$$
$$
\chi\big({}^{\mu}\negmedspace H^0(\phi_f\Fdot[i-1])\big)_\bold p = \chi\big(\phi_f[-1]{}^{\mu}\negmedspace
H^0(\Fdot[i])\big)_\bold p= 
$$
$$
(-1)^{{}^{\dm X}}\big(\operatorname{Ch}({}^{\mu}\negmedspace
H^0(\Fdot[i]))\cdot
\operatorname{im}d\hat f\big)_{(\bold p, d_{\bold p}\hat f)},
$$
where the last equality follows from Theorem 2.1. Now, the result follows from the lemma.\qed

\vskip .5in

\noindent{\it Remark 2.5}. The reader should note that by taking the base ring to be the finite field $\Bbb Z/p$, one
can use Theorem 2.4 to detect $p$-torsion in the integral cohomology of the Milnor fibre (see section 3). Unfortunately, one
cannot quite recover the integral cohomology structure this way; while $\Bbb Z/p$ coefficients would tell one how many direct
summands of the form $\Bbb Z/p^k$ appear in the integral cohomology, they would not distinguish between the various powers
of $p$.

\vskip .2in

It may appear that we have not proved what we claimed we would prove; namely, that we would
conclude Theorem 2.4 under the hypothesis that $\bold p$ was an isolated point in the support $\phi_{f}\Fdot$. However, we
have the following:

\vskip .5in

\noindent{\bf Proposition 2.6}. {\it If $R$ is a p.i.d., then 
$$\big\{\bold p\in X\ |\ f(\bold p)=0, (\bold p, d_\bold p\hat f)\in
|\operatorname{Ch}(\Fdot)|\big\}\ \subseteq\ \operatorname{supp}\phi_f\Fdot.
$$

If $R$ is a field and $\Pdot$ is a perverse sheaf on $X$, then 
$$\big\{\bold p\in X\ |\ f(\bold p)=0, (\bold p, d_\bold p\hat f)\in
|\operatorname{Ch}(\Pdot)|\big\}\ =\ \operatorname{supp}\phi_f\Pdot.
$$
}

\vskip .2in

\noindent{\it Proof}. The first statement is immediate from Theorem 2.10 of [{\bf M1}]. The second statement follows
trivially from Theorem 3.2 of [{\bf M1}].\qed

\vskip .4in

\noindent\S3. {\bf Invisible Strata, Field Coefficients, and the Vanishing Cycles}  

\vskip .2in
We continue with the notation given at the beginning of Section 2, including that $R$ is a p.i.d. 

\vskip .1in

In this section, we begin by discussing how changing the base ring, $R$, affects the vanishing cycles, characteristic cycles, and the
perverse cohomology; this discussion is necessary if we wish to obtain a result which includes possible torsion in the cohomology. This
discussion leads us to define, in 3.3, $\Fdot$-visible strata as those with non-trivial $\Fdot$-Morse modules. We can then prove the main
result of this section, Theorem 3.4, which describes the support of the vanishing cycles; precisely, the result of 3.4 is 
$$
\bigcup_{v\in\Bbb C}\operatorname{supp}\phi_{f-v}\Fdot \ = \ \Big\{\bold x\in X\ |\ (x, d_\bold x\hat
f)\in\bigcup\Sb\Fdot{\text-visible}\\ S_\alpha\endSb\overline{T^*_{{}_{S_\alpha}}\Cal U}\Big\}.
$$

\vskip .3in

\noindent{\bf Basic Results}

\vskip .2in

For each prime ideal $\frak p$ of $R$, let $k_{\frak p}$ denote the field of
fractions of $R/\frak p$, i.e., $k_0$ is the field of fractions of $R$, and for $\frak p\neq 0$, $k_{\frak p} = R/\frak p$.
There are the obvious functors $\delta_{\frak p}: \bold D^b_{{}_\Bbb C}(R_{{}_X})\rightarrow \bold D^b_{{}_\Bbb C}({(k_{\frak
p})}_{{}_X})$, which sends $\Fdot$ to $\Fdot\lotimes (k_{\frak p})^\bullet_{{}_X}$, and $\epsilon_{\frak p}: \bold
D^b_{{}_\Bbb C}({(k_{\frak p})}_{{}_X})\rightarrow\bold D^b_{{}_\Bbb C}(R_{{}_X})$, which considers $k_{\frak p}$-vector spaces as
$R$-modules. 

If $\Adot$ is a
complex of $k_{\frak p}$-vector spaces, we may consider the perverse cohomology of $\Adot$,
${}^{\mu}\negmedspace H^i_{{}_{k_{\frak p}}}(\Adot)$, or the perverse cohomology of $\epsilon(\Adot)$, which we denote by
${}^{\mu}\negmedspace H^i_{{}_R}(\Adot)$. If
$\Adot\in\bold D^b_{{}_\Bbb C}({(k_{\frak p})}_{{}_X})$ and 
$S_{\operatorname{max}}$ is a maximal stratum contained in the support of $\Adot$, then there is a canonical isomorphism
$$\epsilon\big(({}^{\mu}\negmedspace H^i_{{}_{k_{\frak p}}}(\Adot))_{|_{S_\alpha}}\big)\cong ({}^{\mu}\negmedspace
H^i_{{}_{R}}(\Adot))_{|_{S_\alpha}};$$
in particular, $\operatorname{supp}{}^{\mu}\negmedspace H^i_{{}_{k_{\frak p}}}(\Adot) = \operatorname{supp}{}^{\mu}\negmedspace
H^i_{{}_{R}}(\Adot)$.

If $\Fdot\in\bold D^b_{{}_\Bbb C}({R}_{{}_X})$, 
$S_{\operatorname{max}}$ is a maximal stratum contained in the support of $\Fdot$, and $\bold x\in S_{\operatorname{max}}$, 
then for some prime ideal
$\frak p\subset R$ and for some integer $i$, $H^i(\Fdot)_\bold x\otimes k_{\frak p}\neq 0$; it follows that
$S_{\operatorname{max}}$ is also a maximal stratum in the support of $\Fdot\lotimes {(k_{\frak p})}^\bullet_{{}_X}$.
 Thus, 
$$
\operatorname{supp}\Fdot = \bigcup_{\frak p} \operatorname{supp}(\Fdot\lotimes
{(k_{\frak p})}^\bullet_{{}_X})
$$
and so
$$
\operatorname{supp}\Fdot=\bigcup_{i, \frak p}
\operatorname{supp}{}^{\mu}\negmedspace H^i_{{}_{k_{\frak p}}}(\Fdot\lotimes {(k_{\frak p})}^\bullet_{{}_X}),
$$
where the boundedness and constructibility of $\Fdot$ imply that this union is locally finite.

For each prime ideal $\frak p$, there is a natural isomorphism in $\bold D^b_{{}_\Bbb C}(R_{{}_{V(f)}})$ given by
$$
\phi_f\big(\Fdot\lotimes{(k_{\frak p})}^\bullet_{{}_X}\big)\cong\big(\van\big)\lotimes{(k_{\frak p})}^\bullet_{{}_{V(f)}}
$$
(this is a particularly trivial case of the Sebastiani-Thom Isomorphism of [{\bf M3}]), and hence, the stalk cohomology is
given by
$$
H^i\big(\phi_f\big(\Fdot\lotimes{(k_{\frak p})}^\bullet_{{}_X}\big)\big)_\bold x\cong \big(H^i(\van)_\bold x\otimes k_{\frak
p}\big)\ \oplus\ \operatorname{Tor}\big(H^{i+1}(\van)_\bold x,\ k_{\frak p}\big).
$$

\vskip .4in

\noindent{\bf Proposition 3.1}. {\it For all integers $i$ and for all prime ideals $\frak p$ in $R$, the characteristic cycle of the
perverse cohomology of the sheaf of $k_{\frak p}$-vector spaces  $\Fdot[i]\lotimes (k_{\frak p})^\bullet_{{}_X}$ is given by
$$
\operatorname{Ch}({}^{\mu}\negmedspace H^i_{{}_{k_{\frak p}}}(\Fdot\lotimes (k_{\frak p})^\bullet_{{}_X}))=(-1)^{\dm 
X}\sum_\alpha m_{i-d_\alpha}\left[\overline{T^*_{{}_{S_\alpha}}\Cal U}\right],
$$
where $m_{i-d_\alpha}:= {\operatorname{dim}}_{k_{\frak p}}(H^{i-d_\alpha}(\Bbb N_\alpha, \Bbb
L_\alpha; \Fdot)\otimes k_{\frak p}) + {\operatorname{dim}}_{k_{\frak
p}}\operatorname{Tor}(H^{i-d_\alpha+1}(\Bbb N_\alpha, \Bbb
L_\alpha; \Fdot), k_{\frak p})$. }

\vskip .3in

\noindent{\it Proof}. Given the basic facts at the beginning of the section, the proof is exactly that of Lemma 2.3.\qed

\vskip .4in

\noindent{\bf Proposition 3.2} (Continuity of Vanishing Support). {\it Suppose that $\bold x_j\in X$ and $\bold
x_j\rightarrow\bold x$. For each j, let $\hat f_j$ be a complex analytic function locally defined on $\Cal U$ at $\bold
x_j$, such that $\hat f_j(\bold x_j)=0$, and let $f_j:=\hat {f_j}_{|_X}$. Suppose that $\bold
x_j\in\operatorname{supp}\phi_{f_j}\Fdot$ and that $d_{\bold x_j}\hat f_j\rightarrow d_\bold x\hat f$.

\vskip .1in

Then, the following three equivalent conclusions hold:

\vskip .1in

\noindent a)\hskip .2in $\bold x\in\operatorname{supp}\phi_{f-f(\bold x)}\Fdot$;

\vskip .1in

\noindent b)\hskip .2inthere exists an integer
$i$ and the prime ideal
$\frak p$ in $R$ such that
$$
\bold x\in\operatorname{supp}\phi_{f-f(\bold x)}{}^{\mu}\negmedspace H^i_{{}_{k_{\frak p}}}(\Fdot\lotimes (k_{\frak
p})^\bullet_{{}_X}); 
$$
and
\vskip .1in

\noindent c)\hskip .2in there exists an integer
$i$ and the prime ideal
$\frak p$ in $R$ such that
$$
(\bold x, d_\bold x\hat f)\in\big|\operatorname{Ch}({}^{\mu}\negmedspace H^i_{{}_{k_{\frak p}}}(\Fdot\lotimes (k_{\frak
p})^\bullet_{{}_X}))\big|.
$$ }

\vskip .3in

\noindent{\it Proof}. Let us prove the equivalence of a), b), and c) first. By 2.6, b) and c) are equivalent. Now, a) and b) are
equivalent because
$$
\operatorname{supp}\phi_{f-f(\bold x)}\Fdot = \bigcup_{i, \frak p}\operatorname{supp}{}^{\mu}\negmedspace H^i_{{}_{k_{\frak
p}}}(\phi_{f-f(\bold x)}\Fdot\lotimes {(k_{\frak p})}^\bullet_{{}_{V(f-f(\bold x))}}) =$$ 
$$\bigcup_{i, \frak
p}\operatorname{supp}{}^{\mu}\negmedspace H^i_{{}_{k_{\frak p}}}\big(\phi_{f-f(\bold x)}(\Fdot\lotimes {(k_{\frak
p})}^\bullet_{{}_{X}})\big)=\bigcup_{i, \frak
p}\operatorname{supp}\phi_{f-f(\bold x)}{}^{\mu}\negmedspace H^i_{{}_{k_{\frak p}}}\big(\Fdot\lotimes {(k_{\frak
p})}^\bullet_{{}_{X}}\big).
$$

\vskip .2in

We will now prove that the assumptions of the first paragraph imply c).
As above, we have
$$\operatorname{supp}\phi_{f_j}\Fdot = \bigcup_{i, \frak
p}\operatorname{supp}\phi_{f_j}{}^{\mu}\negmedspace H^i_{{}_{k_{\frak p}}}\big(\Fdot\lotimes {(k_{\frak
p})}^\bullet_{{}_{X}}\big).
$$

Now, near $\bold x$, $\operatorname{supp}\Fdot$ is a finite union $\bigcup_{i, \frak p}
\operatorname{supp}{}^{\mu}\negmedspace H^i_{{}_{k_{\frak p}}}(\Fdot\lotimes {(k_{\frak p})}^\bullet_{{}_X})$; thus, by
taking a subsequence, we may assume that there is one $i$ and one $\frak p$ such that, for all $j$, 
$\bold x_j\in \operatorname{supp}\phi_{f_j}{}^{\mu}\negmedspace H^i_{{}_{k_{\frak p}}}\big(\Fdot\lotimes {(k_{\frak
p})}^\bullet_{{}_{X}}\big)$. It follows from 2.6 that $(\bold x_j, d_{\bold x_j}\hat f_j)\in
\big|\operatorname{Ch}\big({}^{\mu}\negmedspace H^i_{{}_{k_{\frak p}}}\big(\Fdot\lotimes {(k_{\frak
p})}^\bullet_{{}_{X}}\big)\big)\big|$. As $\big|\operatorname{Ch}\big({}^{\mu}\negmedspace H^i_{{}_{k_{\frak p}}}\big(\Fdot\lotimes {(k_{\frak
p})}^\bullet_{{}_{X}}\big)\big)\big|$ is closed, we conclude that $(\bold x, d_\bold x\hat f)\in \big|\operatorname{Ch}\big({}^{\mu}\negmedspace H^i_{{}_{k_{\frak p}}}\big(\Fdot\lotimes {(k_{\frak
p})}^\bullet_{{}_{X}}\big)\big)\big|$.\qed

\vskip .4in

Looking at Theorem 2.4, and Propositions 3.1 and 3.2, we see that any stratum $S_\alpha$ for which $H^*(\Bbb N_\alpha, \Bbb
L_\alpha;\Fdot)=0$ is essentially irrelevant as far as vanishing cycles are concerned. Hence, we make the following definition.

\vskip .3in

\noindent{\bf Definition 3.3}. A stratum $S_\alpha$ is {\it $\Fdot$-invisible} if $H^*(\Bbb N_\alpha, \Bbb
L_\alpha;\Fdot)=0$. Otherwise, we say that $S_\alpha$ is {\it $\Fdot$-visible}.

\vskip .3in

The basic principle is that $\Fdot$-invisible strata do not contribute any cohomology in Morse Theory arguments. As an example, we have 

\vskip .3in

\noindent{\bf Theorem 3.4}. {\it
$$
\bigcup_{v\in\Bbb C}\operatorname{supp}\phi_{f-v}\Fdot \ = \ \Big\{\bold x\in X\ |\ (x, d_\bold x\hat
f)\in\bigcup\Sb\Fdot{\text-visible}\\ S_\alpha\endSb\overline{T^*_{{}_{S_\alpha}}\Cal U}\Big\}.
$$

}

\vskip .3in

\noindent{\it Proof}. This is immediate from the equivalence of a) and c) in 3.2, and the description of the characteristic cycle
given in 3.1.\qed

\vskip .4in

\noindent{\it Remark 3.5}. The union on the left side above is not just locally finite, but, in fact, locally consists of a single
support, i.e., near a point
$\bold p\in X$,
$\operatorname{supp}\phi_{f-v}\Fdot =\emptyset$ unless $v=f(\bold p)$.

We should also point out an elementary, but important, relation with the real structure of the conormal spaces. In the situation of
Theorem 3.4, $d_\bold x\hat f\in \big(\overline{T^*_{{}_{S_\alpha}}\Cal U}\big)_\bold x$ if and only if $d_\bold x(\operatorname{Re}\hat
f)\in
\big(\overline{T^*_{{}_{S_\alpha}}\Cal U}\big)_\bold x$ (considered with its real structure); this follows from the Cauchy-Riemann
equations. Of course,
$d_\bold x\hat f\in
\big(\overline{T^*_{{}_{S_\alpha}}\Cal U}\big)_\bold x$ if and only if there exists $a+bi\neq 0$ such that $d_\bold x[(a+bi)\hat
f]=(a+bi)d_\bold x\hat f\in
\big(\overline{T^*_{{}_{S_\alpha}}\Cal U}\big)_\bold x$. It follows that $d_\bold x\hat f\not\in
\big(\overline{T^*_{{}_{S_\alpha}}\Cal U}\big)_\bold x$ if and only if, for all $(a,b)\in\Bbb R^2-\{\bold 0\}$, $a\,d_\bold
x(\operatorname{Re}\hat f)+b\,d_\bold x(\operatorname{Im}\hat f)\not\in \big(\overline{T^*_{{}_{S_\alpha}}\Cal U}\big)_\bold x$. This will
be importance to us in section 5.

\vskip .4in

\noindent\S4. {\bf Microlocal Basics}  

\vskip .2in

In this section, we will discuss conormal geometry, and provide a down-to-Earth discussion of a number of microlocal
results of Kashiwara and Schapira. We continue with our earlier notations, including that the base ring $R$ is a p.i.d.

\vskip .2in

We will always work inside an open subset $\Cal U$ in $\Bbb R^{m+1}$. Let $T^*\Cal U@>\pi>>\Cal U$ denote the conormal
bundle. This is isomorphic to $\Cal U\times\Bbb R^{m+1}@>\operatorname{pr}>>\Cal U$; we use $x_0, \dots, x_m, \xi_0, \dots,
\xi_m$ for coordinates on $T^*\Cal U$. We use $\Bbb R\Bbb P(T^*\Cal U)$ to denote the total space of the projectivized cotangent bundle;
hence, $\Bbb R\Bbb P(T^*\Cal U)\cong \Cal U\times \Bbb R\Bbb P^m$.

If $M$ is a smooth submanifold of $\Cal U$, then the {\it conormal bundle to $M$ in $\Cal U$} is given by 
$$T^*_{{}_M}\Cal U
:=\{(\bold p,
\eta)\in T^*\Cal U\ |\
\bold p\in M,\ 
\eta(T_\bold p M)=0\}.$$
Using this notation, the zero-section is given by $T^*_{{}_\Cal U}\Cal U$. In addition, two submanifolds $M,
N\subseteq\Cal U$ intersect transversely if and only if $T^*_{{}_M}\Cal U \cap T^*_{{}_N}\Cal U \subseteq T^*_{{}_\Cal
U}\Cal U$.

There is a canonical $1$-form, $\alpha$, on $T^*\Cal U$, which we will describe at each point. Let $(\bold p, \eta)\in
T^*\Cal U$. We want to describe the linear function
$\alpha_{(\bold p, \eta)}: T_{(\bold p, \eta)}(T^*\Cal U)\rightarrow\Bbb R$. Note that $\eta: T_{\bold p}\Cal
U\rightarrow\Bbb R$ and that $d_{(\bold p, \eta)}\pi: T_{(\bold p, \eta)}(T^*\Cal U)\rightarrow T_{\bold p}\Cal U$. Thus,
we may define $\alpha_{(\bold p, \eta)}:= \eta\circ d_{(\bold p, \eta)}\pi$. One easily verifies that, in coordinates,
$\alpha$ is given by $\xi_0dx_0 + \dots + \xi_mdx_m$.

Let $\Bbb R^+$ denote the strictly positive real numbers. A subset $S$ of $T^*\Cal U$ is {\it $\Bbb R^+$-conic} provided
that, for all $\bold p\in\Cal U$, for all $r\in \Bbb R^+$ and $\eta\in \pi^{-1}(\bold p)$, $r\eta\in S$, i.e., the fibres
of
$S$ are closed under positive scalar multiplication. While the zero vector need not be in a fibre if a conic subset, if
$S$ is closed and conic, then the zero vector is in every non-empty fibre. In addition, it is easy to prove that the image
under $\pi$ of a closed, conic subset of $T^*\Cal U$ is a closed subset of $\Cal U$. If $S$ is $\Bbb R^+$-conic, then we may consider its
projectivization $\Bbb R\Bbb P(S)\subseteq \Bbb R\Bbb P(T^*\Cal U)$.

\vskip .2in

We need to define and describe isotropic subsets of $T^*\Cal U$.

\vskip .2in

\noindent{\bf Definition 4.1}. A subset $S\subseteq T^*\Cal U$ is {\it isotropic} if $\alpha_{|_S}=0$.

\vskip .4in

The following proposition provides a nice characterization of isotropic sets; it is 8.3.10 of [{\bf K-S3}].

\vskip .2in

\noindent{\bf Proposition 4.2}. {\it Let $S$ be an $\Bbb R^+$-conic subanalytic subset of $T^*\Cal U$. Then, the
following are equivalent:

\vskip .1in

\noindent i)  \ \ $S$ is isotropic;

\vskip .1in

\noindent ii) \ there exists a locally finite family $\{W^j\}$ of subanalytic subsets of $\Cal U$ such that $S\subseteq
\bigcup_j
\overline{T^*_{{}_{W^j_{\operatorname{reg}}}}\Cal U}$;

\vskip .1in

\noindent iii) there exists a finite family $\{W^j\}$ of subanalytic submanifolds of $\Cal U$ such that $W^j\subseteq
\pi(S)$ and 
$S\subseteq
\bigcup_j
T^*_{{}_{W^j}}\Cal U$.
}

\vskip .4in

Our primary interest in isotropic subsets stems from the following theorem ([{\bf K-S3}], 8.3.12), which
essentially says that if $S$ is a closed, $\Bbb R^+$-conic, subanalytic, isotropic subset of $T^*\Cal U$, then the
critical values of any proper real analytic function, relative to $S$, are discrete.

\vskip .2in

\noindent{\bf Theorem 4.3} (Microlocal Bertini-Sard Theorem). {\it Let $\phi:\Cal U\rightarrow\Bbb R$ be a real
analytic function and let $S\subseteq T^*\Cal U$ be a closed, $\Bbb R^+$-conic, subanalytic, isotropic subset. Assume
that $\phi$ is proper on $\pi(S)$.

Then, the set $\{t\in\Bbb R\ |\ \text{there exists } \bold x\in\Cal U\text{ such that } t=\phi(\bold x)\text{ and
}(\bold x, d_{\bold x}\phi)\in S\}$ is discrete.
}

\vskip .5in

If $E\subseteq T^*(\Cal U)$, then $E^a$ denotes the image of $E$ by the antipodal map, i.e., 
$$E^a:= \{(\bold p, -\eta)\
|\ (\bold p, \eta)\in E\}.$$

\vskip .2in

If $A$ and $B$ are $\Bbb R^+$-conic subsets of $T^*\Cal U$, then $A+B$ is the $\Bbb R^+$-conic subset of $T^*\Cal U$ which lies
over
$\pi(A)\cap \pi(B)$ and such that
$$
\pi^{-1}(\bold x)\cap (A+B) = \{(\bold x, a+b)\ |\ (\bold x, a)\in A, \ (\bold x, b)\in B\}.
$$
If $A$ and $B$ are closed, then certainly the fibres of $A+B$ are closed; however, $A+B$, itself, need not be closed (but see 4.6 below).

\vskip .4in

We need to describe the $\widehat+$ operation of [{\bf K-S3}], 6.2.3.v. While the definition of $\widehat+$ is somewhat
difficult to unravel, the characterization of the operation given in 6.2.8.ii of [{\bf K-S3}] is easy to understand; hence, we will use
this as our definition.

\vskip .2in

\noindent{\bf Definition 4.4}. Let $A$ and $B$ be two closed $\Bbb R^+$-conic subsets of $T^*\Cal U$. Then,
$A\widehat+B$ is the subset of $T^*\Cal U$ defined by: $(\bold p, \eta)\in A\widehat+B$ if and only if there exist
sequences $(\bold x_n, \sigma_n)\in A$ and $(\bold y_n, \tau_n)\in B$ such that $\bold x_n\rightarrow\bold p$, $\bold
y_n\rightarrow\bold p$, $\sigma_n+\tau_n\rightarrow \eta$, and $|\bold x_n-\bold y_n|\cdot|\sigma_n|\rightarrow 0$.

\vskip .4in

\noindent{\it Remark 4.5}. It follows from the definition that $A\widehat+B$ is closed, and lies over the closed set
$\pi(A)\cap\pi(B)$. Also, while the definition may appear to be asymmetric, in fact, the conditions imply that $|\bold
x_n-\bold y_n|\cdot|\tau_n|\rightarrow 0$. Note that $A+B\subseteq A\widehat+B$.

\vskip .4in

The following is a combination of Lemma 5.4.7 and Corollary 8.3.18.i of [{\bf K-S3}].

\vskip .2in

\noindent{\bf Proposition 4.6}. Let $A$ and $B$ be two closed, $\Bbb R^+$-conic subsets of
$T^*\Cal U$. 

\vskip .1in
 \noindent i)\hskip .2in Then, $A\widehat+B$ is also closed and $\Bbb R^+$-conic.

\vskip .1in

 \noindent ii)\hskip .18in If $A\cap B^a\subseteq T_{{}_{\Cal U}}^*\Cal U$, then $A+B$ 
is closed and $\Bbb R^+$-conic.

\vskip .1in

 \noindent iii)\hskip .16in If $A$ and $B$ are  subanalytic, isotropic subsets of
$T^*\Cal U$, then $A\widehat+B$ is a closed, $\Bbb R^+$-conic, subana-

\noindent \hbox{}\hskip .33in lytic, isotropic subset of
$T^*\Cal U$.

\vskip .4in

For $\Fdot\in D^b(\Cal U)$, we need to define the micro-support, $SS(\Fdot)$, of $\Fdot$; the micro-support consists of
a set of covectors which point in directions in which $\Fdot$ ``does not propagate''-- that is, directions in which the
hypercohomology of $\Fdot$ changes locally.

\vskip .2in

\noindent{\bf Definition 4.7}. The {\it micro-support of $\Fdot$}, $SS(\Fdot)$, is the subset of $T^*\Cal U$ defined by
the following: $(\bold p, \eta)\not\in SS(\Fdot)$ if and only if there exists an open neighborhood $\Omega$ of $(\bold
p, \eta)$ in $T^*\Cal U$ such that for all $\bold x\in \Cal U$, for all real $C^1$ functions $\psi$ defined in a
neighborhood of $\bold x$ such that $(\bold x, d_{\bold x}\psi)\in\Omega$ and $\psi(\bold x)=0$,
$$
\big(R\Gamma_{{}_{\{\bold y\in\Cal U|\psi(\bold y)\geqslant 0\}}}(\Fdot)\big)_\bold x = 0.
$$

\vskip .4in

The following
proposition is 5.1.3.i of [{\bf K-S3}], combined with part of 8.4.2.

\vskip .2in

\noindent{\bf Proposition 4.8}. {\it Let $\Fdot\in D^b(\Cal U)$. Then, $SS(\Fdot)$ is a closed, $\Bbb R^+$-conic subset
of $T^*\Cal U$ and \hbox{$SS(\Fdot)\cap T^*_{{}_{\Cal U}}\Cal U= \operatorname{supp}\Fdot\times\{0\}$}.

In addition, if $\Fdot$ is $\Bbb R$-constructible, then $SS(\Fdot)$ is a subanalytic, isotropic subset in $T^*\Cal U$.
}

\vskip .4in

The following lemma is 5.4.19 of [{\bf K-S3}]. It tells us that if a function produces no infinitesimal changes in
$\Fdot$, then it produces no global changes in $\Fdot$.

\vskip .2in

\noindent{\bf Lemma 4.9}. {\it Let $\Fdot\in D^b(\Cal U)$ and let $\phi:\Cal U\rightarrow\Bbb
R$ be a $C^1$ function such that $\phi_{|_{\operatorname{supp}(\Fdot)}}$ is proper. Let $a, b\in\Bbb R$ with $a<b$.

If, for all $\bold x\in\phi^{-1}([a,b))$,  $(\bold x, d_\bold x\phi)\not\in SS(\Fdot)$, then the natural morphisms
$$
R\Gamma\big(\phi^{-1}((-\infty,b); \Fdot)\big)\rightarrow R\Gamma\big(\phi^{-1}((-\infty,a]; \Fdot)\big)\rightarrow
R\Gamma\big(\phi^{-1}((-\infty,a); \Fdot)\big)
$$
are isomorphisms.
}

\vskip .4in

The following is part of Lemma 8.4.7 of [{\bf K-S3}].

\vskip .3in

\noindent{\bf Lemma 4.10}. {\it Let $X$ be a real analytic space, $\Fdot\in D^b_{{}_{\Bbb R}}(X)$, and $r:X\rightarrow\Bbb
R^n$ a real analytic function. Assume that $r_{|_{\operatorname{supp}\Fdot}}$ is proper. Then, for all
sufficiently small $\epsilon>0$, the inclusion maps induce natural isomorphisms
$$
R\Gamma(r^{-1}(B_\epsilon);\ \Fdot)\ \cong\ R\Gamma(r^{-1}\big({\overset\circ\to B}_\epsilon\big);\ \Fdot)\ \cong\
R\Gamma(r^{-1}(\bold 0);\ \Fdot).
$$
}

\vskip .4in

Our application of 4.9 will be to the case where we begin with a complex $\Fdot\in D^b_{{}_\Bbb C}(\Cal U)$, then
restrict $\Fdot$ to a closed ball, and then push it forward, back onto $\Cal U$. Thus, we need to know how restricting
to a closed ball, and then pushing forward, affects the micro-support. For this, we need two or three more definitions.

\vskip .2in

\noindent{\bf Definition 4.11}. Suppose that $Z\subseteq\Cal U$ is a closed submanifold with boundary, of dimension
equal to that of $\Cal U$. For each point $\bold x\in\partial Z$, there is a unique, inwardly-pointing unit vector,
$\bold n_\bold x$, which is normal to $\partial Z$ at $\bold x$. The standard inner-product (dot product) gives a
well-defined linear form $\eta_\bold x(\bold v):= \bold v\cdot \bold n_\bold x$. 

The subset $N^*(Z)\subseteq T^*(\Cal U)$ is defined
by  $\pi(N^*(Z)) = \Cal U$ and, for all $\bold x\in\Cal U$,
$$
\pi^{-1}(\bold x)\cap N^*(Z) = \cases \{(\bold x, 0)\}, & \text{ if } \bold x\not\in\partial Z\\\{(\bold x, r\eta_\bold
x)\ |\ r\in \{0\}\cup\Bbb R^+\},& \text{ if }
\bold x\in\partial Z.
\endcases
$$
(See 5.3.6 of [{\bf K-S3}].) Hence, $N^*(Z)$ consists of the zero-section of $T^*\Cal U$, together with the
``inwardly-pointing conormals'' along the boundary of $Z$. Note that $N^*(Z)$ is closed and $\Bbb R^+$-conic. In
addition, if $Z$ and $\partial Z$ are subanalytic, then Proposition 4.2 implies that $N^*(Z)$ is isotropic.

\vskip .5in

Recall that if $j:Z\hookrightarrow\Cal U$ is a closed inclusion, then $\Fdot_Z\cong j_!j^*\Fdot\cong j_*j^*\Fdot$;
therefore, Proposition 5.4.8.b.ii of [{\bf K-S3}] yields the following:

\vskip .2in

\noindent{\bf Proposition 4.12}. {\it Let $\Fdot\in D^b(\Cal U)$. Let $j:Z\hookrightarrow\Cal U$ be a closed submanifold
with boundary, of dimension equal to that of $\Cal U$.

If $SS(\Fdot)\cap N^*(Z)^a\subseteq T^*_{{}_{\Cal U}}\Cal U$, then $SS(j_*j^*\Fdot)\subseteq N^*(Z) + SS(\Fdot)$.
}

\vskip .2in

One should think of the condition that $SS(\Fdot)\cap N^*(Z)^a\subseteq T^*_{{}_{\Cal U}}\Cal U$ as a transversality
condition between $\partial Z$ and some ``virtual strata'' whose conormals form $SS(\Fdot)$.

\vskip .4in

In the constructible complex analytic setting, we are able to link the micro-support with the visible strata of Definition 3.3.

\vskip .4in

\noindent{\bf Theorem 4.13}. {\it Suppose that $\Fdot\in D^b_{{}_\Bbb C}(\Cal U)$ and that $\{S_\alpha\}$ is a complex analytic Whitney
stratification of $\Cal U$ with connected strata. Then,
$$
SS(\Fdot)= \bigcup\Sb\Fdot{\text-visible}\\ S_\alpha\endSb\overline{T^*_{{}_{S_\alpha}}\Cal U}.
$$

}

\vskip .3in

\noindent{\it Proof}. Let $(\bold x, \eta)\in T^*\Cal U$. By Proposition 3.1, $$\bigcup\Sb\Fdot{\text-visible}\\
S_\alpha\endSb\overline{T^*_{{}_{S_\alpha}}\Cal U} \ =\ \bigcup_{i,\frak p}\big|\operatorname{Ch}({}^{\mu}\negmedspace H^i_{{}_{k_{\frak
p}}}(\Fdot\lotimes (k_{\frak p})^\bullet_{{}_X}))\big|,\tag{$\dagger$}
$$
where the union is over all integers
$i$ and the prime ideals
$\frak p$ in $R$.

 By 8.6.4 of [{\bf K-S3}],
$(\bold x,
\eta)\in SS(\Fdot)$ if and only if there exist
$\bold x_j\in \Cal U$ and locally defined complex analytic $\hat f_j$  such that $f_j(\bold x_j)=0$, $H^*(\phi_{f_j}\Fdot)_{\bold x_j}\neq
0$, and  $d_{\bold x_j}\hat f_j\rightarrow \eta$; clearly, this is equivalent to: there exist
$\bold x_j\in \Cal U$ and locally defined complex analytic $\hat f_j$  such that $f_j(\bold x_j)=0$, $\bold
x_j\in\operatorname{supp}\phi_{f_j}\Fdot$, and  $d_{\bold x_j}\hat f_j\rightarrow \eta$. By the Continuity of Vanishing Support (3.2),
this implies that there exists an integer
$i$ and the prime ideal
$\frak p$ in $R$ such that
$$
(\bold x, \eta)\in\big|\operatorname{Ch}({}^{\mu}\negmedspace H^i_{{}_{k_{\frak p}}}(\Fdot\lotimes (k_{\frak
p})^\bullet_{{}_X}))\big|.
$$
Hence, $\dsize SS(\Fdot)\subseteq \bigcup\Sb\Fdot{\text-visible}\\ S_\alpha\endSb\overline{T^*_{{}_{S_\alpha}}\Cal U}$.

Conversely, if $\dsize (\bold x, \eta)\in\bigcup\Sb\Fdot{\text-visible}\\ S_\alpha\endSb\overline{T^*_{{}_{S_\alpha}}\Cal U}$, then, by
$(\dagger)$ and 3.2, $\bold x\in\operatorname{supp}\phi_{{}_{L-L(\bold x)}}\Fdot$ for any linear form $L$ such that $d_\bold x L=\eta$.
Now, 8.6.4 of [{\bf K-S3}] immediately implies that $(\bold x, \eta)\in SS(\Fdot)$.\qed

\vskip .4in

\noindent{\it Remark 4.14}. Theorem 4.13 immediately implies a stronger version of itself. One does not need to begin with a Whitney
stratification, but merely any stratification for which the normal data of strata with respect to $\Fdot$ is ``well-defined'', i.e.,
stratifications in which the normal data in normal slices to strata locally trivializes along the strata. Such stratifications only
require refinement by including
$\Fdot$-invisible strata in order to obtain a Whitney stratification. 

This is essentially what is required by Brian\c con, Maisonobe, and Merle in [{\bf BMM}], where they use stratifications which satisfy
Whitney's condition a) and the {\it property of local stratified triviality}.

\vskip .4in

\noindent{\bf Corollary 4.15}. {\it In the notation and situation of Theorem 3.4,
$$
\bigcup_{v\in\Bbb C}\operatorname{supp}\phi_{f-v}\Fdot \ = \ \Big\{\bold x\in X\ |\ (\bold x, d_\bold x\hat
f)\in SS(\Fdot)\Big\}.
$$
}

\vskip .3in

\noindent{\it Proof}. This follows immediately from 3.4 and 4.13, once we extend $\Fdot$ by zero to all of $\Cal U$.\qed

\vskip .4in

\noindent\S5. {\bf The Microlocal Proof}  

\vskip .2in

We now return to the setting of Theorem 1.1. Let $X$ be a complex analytic space
embedded in an open subset $\Cal U$ of $\Bbb C^{n+1}$, let $\Cal S:=\{S_\alpha\}$ be a complex Whitney stratification of $X$ with
connected strata, let
$\Fdot\in D^b_{{}_\Cal S}(X)$, where the base ring is a p.i.d., and let $\hat f:\Cal U\rightarrow\Bbb C$ be complex analytic. Let $f:=
\hat f_{|_{X}}$. Suppose that $\bold 0\in X$ and $f(\bold 0)=0$.

If we let $j:X\hookrightarrow\Cal U$ denote the inclusion, then the problem of studying $\Fdot$, $f$, $X$, and $\Cal S$ is equivalent to
studying $j_!\Fdot$, $\hat f$, $\Cal U$, and $\Cal S\cup\{\Cal U-X\}$. {\bf Therefore, throughout this section, we will assume without
loss of generality, that
$X=\Cal U$, and so we write $f$ in place of $\hat f$.}

\vskip .3in

We want to know how to alter the proof of Theorem 1.1 when we weaken our hypothesis, and only suppose
that $\bold 0$ is an isolated point in the support of $\phi_f\Fdot$, rather than being an isolated stratified critical point. 

In fact, all we need are two lemmas: one which says the that a general complex linear perturbation of $f$ has isolated stratified
critical points, and one which says that any possible stratified critical points on the boundary of a small ball are irrelevant.

\vskip .4in

\noindent{\bf Lemma 5.1}. {\it Let $L:X\rightarrow\Bbb C$  be a generic linear form. Then,
there exists $\epsilon>0$ such that, for all sufficiently small $\delta>0$, for all $t\in\Bbb D^*_\delta$, $f+tL$ has isolated
stratified critical points inside
$\overset\circ\to B_\epsilon$ and no critical points which are contained in $\partial B_\epsilon$.
}

\vskip .3in

\noindent{\it Proof}. We believe that this is well-known -- even in our case, where $f$ may have stratified non-isolated critical points;
however, for lack of a convenient reference in this generality, we supply a proof. 

Fix a stratum $S_\alpha\neq \{\bold 0\}$. Let $L$ be generic enough so that $d_\bold 0 L\not\in \big(\overline{T^*_{{}_{S_\alpha}}X}\big)_\bold 0$.
Also, choose $L$ so generic that the relative polar curve,
$\Gamma^1_{f_{|_{{\overline{S}}_\alpha}}, L}$, is purely $1$-dimensional at the origin, and so that $L$ is finite near the origin when
restricted to the relative polar curve (that this is possible is well-known and is proved in many places; see, for instance, [{\bf L-T}],
4.2.1 and [{\bf Te}], 4.1.3.2.).

For any complex number $t\neq 0$, a critical point of $f+tL$ on $S_\alpha$ is either a critical point of $f_{|_{S_\alpha}}$ or is a point
on the relative polar curve. If we had an entire curve, containing the origin in its closure, of critical points of both 
$(f+tL)_{|_{S_\alpha}}$ and $f_{|_{S_\alpha}}$, then we would have $d_\bold 0 L\in \big(\overline{T^*_{{}_{S_\alpha}}X}\big)_\bold 0$ --
a contradiction. Thus, any curve of critical points of $(f+tL)_{|_{S_\alpha}}$ at the origin must be a component of
$\Gamma^1_{f_{|_{{\overline{S}}_\alpha}}, L}$. We will show that this can only happen for a finite number of ``bad'' $t$.

Let $\omega(u)$ be an analytic paramerization of a component of $\Gamma^1_{f_{|_{{\overline{S}}_\alpha}}, L}$ such that $\omega(0)=\bold
0$ and, for $u\neq 0$, $\omega(u)\in\Sigma(f+tL)_{|_{S_\alpha}}$. Then, $d_{\omega(u)}(f+tL)(\omega^\prime(u))\equiv 0$ for
$u\neq 0$. Therefore, $f(\omega(u))+tL(\omega(u))\equiv 0$, and so, as $L(\omega(u))\not\equiv 0$ by our finiteness assumption, it
follows that 
$\dsize\lim_{u\rightarrow 0}\frac{-f(\omega(u))}{L(\omega(u))}$ exists and is equal to $t$. Hence, there are at most as many ``bad'' $t$
values as there are components of $\Gamma^1_{f_{|_{{\overline{S}}_\alpha}}, L}$. The desired conclusion follows.\qed

\vskip .4in

In the following lemma, we write $T^*_{{}_{V(f)}}X$ in place of $\overline{T^*_{{}_{(V(f))_{\operatorname{reg}}}}X}$. We also write
$\Sigma f$ for the critical locus of $f$, which is unambiguous since we are assuming that the domain of $f$ is an open
subset of affine space. Finally, in some places below, we will identify $T^*X$ with $X\times\Bbb C^{n+1}\cong X\times\Bbb R^{2n+2}$, so
that we may use the real metric and  minors of matrices.

\vskip .4in

\noindent{\bf Lemma 5.2}. {\it Suppose
that $\bold 0$ is an isolated point in the support of $\phi_f\Fdot$. For complex numbers $t$ and $v$, and linear forms $L$, let
$\Psi_{L,t,v}:X\rightarrow\Bbb R$ be given by
$\Psi_{L,t,v}(\bold z):=|f(\bold z)+tL(\bold z)-v|^2$. 

For all linear forms $L$, for all sufficiently small
$\epsilon>0$, if
$j:B_\epsilon\hookrightarrow X$ denotes the inclusion,  then there exists an open neighborhood $\Cal W$ of $\partial B_\epsilon\cap
V(f)$ in $\partial B_\epsilon$ and an open neighborhood $\Cal Y$ of $0$ in $\Bbb C$ such that, for all
$(\bold x, t, v)\in \Cal W\times\Cal Y\times\Bbb C$, $(\bold x, d_\bold x \Psi_{L,t,v})\not\in
SS(j_*j^*\Fdot)$ unless $\Psi_{L,t,v}(\bold x)= 0$.
}

\vskip .3in

\noindent{\it Proof}. Fix a linear form $L$. We will make a continuity argument to show that there exist such $\Cal W$ and $\Cal Y$ such that, for
all
$(\bold x, t)\in \Cal W\times\Cal Y$, for all $(a,b)\in\Bbb R^2-\{\bold 0\}$,
$$\big(\bold x, \ a\,d_\bold x(\operatorname{Re} (f+tL))+b\,d_\bold x(\operatorname{Im} (f+tL))\big)\not\in SS(j_*j^*\Fdot).$$
This yields the desired result, since $$d_\bold x\Psi_{L,t,v} \ =\ 2\operatorname{Re} (f(\bold x)+tL(\bold x)-v)\,d_\bold x(\operatorname{Re}
(f+tL))+2\operatorname{Im} (f(\bold x)+tL(\bold x)-v)\,d_\bold x(\operatorname{Im} (f+tL)).$$

\vskip .3in

By Corollary 4.15, there exists $\epsilon_1>0$ such that, for all $\bold x\in B_{\epsilon_1}-\{\bold 0\}$, $(\bold x, d_\bold x f)\not\in
SS(\Fdot)$, i.e.,  for all $\bold x\in (B_{\epsilon_1}-\{\bold 0\})\cap \operatorname{supp}\Fdot$, $d_\bold x f\not\in \big(SS(\Fdot)\big)_\bold
x$; in particular, for all $\bold x\in (B_{\epsilon_1}-\{\bold 0\})\cap \operatorname{supp}\Fdot$, $d_\bold x f\neq 0$. By Remark 3.5, we
conclude that, for all $\bold x\in (B_{\epsilon_1}-\{\bold 0\})\cap \operatorname{supp}\Fdot$, for all $(a, b)\in\Bbb R^2-\{\bold 0\}$,
$a\,d_\bold x(\operatorname{Re} f)+b\,d_\bold x(\operatorname{Im} f)\not\in \big(SS(\Fdot)\big)_\bold x$.

\vskip .2in

Now consider $SS(\Fdot)\ \widehat +\ T^*_{{}_{V(f)}}X$; this is a closed, $\Bbb R^+$-conic (actually, $\Bbb C$-conic), isotropic, subanalytic
subset of $T^*X$. Let $r:X\rightarrow\Bbb R$ be given by $r(\bold z):= |\bold z|^2$. Then, the microlocal
Bertini-Sard Theorem (4.3) implies that there exists $\epsilon_2>0$ such that, for all $\bold x\in r^{-1}(0,\epsilon_2^2)$, $$(\bold x,
d_\bold x r)\not\in SS(\Fdot)\ \widehat +\ T^*_{{}_{V(f)}}X\ \supseteq \ SS(\Fdot)\ +\ T^*_{{}_{V(f)}}X.\tag{$\dagger$}$$

\vskip .2in

Let $\epsilon_0=\operatorname{min}\{\epsilon_1, \epsilon_2\}$. Let $U(SS(\Fdot))$ be the covectors of $SS(\Fdot)$ of unit length, i.e., 
$U(SS(\Fdot)):=
\{(\bold x,\omega)\in SS(\Fdot)\ |\
\ ||\omega||=1\}$ (we have used the identification $T^*X\cong X\times\Bbb R^{2n+2})$. From the previous two paragraphs, we conclude that, for all
$\bold x\in (B_{\epsilon_0}-\{\bold 0\})\cap \operatorname{supp}\Fdot$, for all $\omega \in \big(U(SS(\Fdot))\big)_\bold x$, the covectors 
$d_\bold x(\operatorname{Re} f)$, $d_\bold x(\operatorname{Im} f)$, $d_\bold x r$, and $\omega$ are linearly independent over $\Bbb R$.

\vskip .3in

Consider the function $M: U(SS(\Fdot))\times\Bbb C\rightarrow X\times \Bbb R^k$ (for an appropriate $k$) given by 
$$M(\bold x,\omega, t):=\left(\bold
x, \ 4\times 4 \text{ minors of } \pmatrix d_\bold x(\operatorname{Re} (f+tL))\\d_\bold x(\operatorname{Im} (f+tL))\\d_\bold x
r\\ \omega\endpmatrix\right).
$$
Then, $M$ is a continuous function, and, if $0<\epsilon<\epsilon_0$,
$$
M^{-1}(X\times\{\bold 0\})\ \cap\ \Big(\big[\big((\partial B_\epsilon\cap V(f))\times\Bbb R^{2n+2}\big)\cap U(SS(\Fdot))\big] \times\{0\}\Big)\
=\ \emptyset.
$$
By normality, it follows that there exists an open neighborhood $\Cal Z$ of $\big[\big((\partial B_\epsilon\cap V(f))\times\Bbb R^{2n+2}\big)\cap
U(SS(\Fdot))\big] \times\{0\}$ in $U(SS(\Fdot))\times\Bbb C$ such that, for all $(\bold x, \omega, t)\in \Cal Z$, $d_\bold x(\operatorname{Re}
(f+tL))$, $d_\bold x(\operatorname{Im} (f+tL))$, $d_\bold x r$, and $\omega$ are linearly independent over $\Bbb R$. 

Using that $SS(\Fdot)$ is
$\Bbb R^+$-conic and that $\partial B_\epsilon\cap
V(f)\cap\operatorname{supp}\Fdot$ is compact, we immediately conclude that there exist an open neighborhood $\Cal W^\prime$ of $\partial
B_\epsilon\cap V(f)\cap\operatorname{supp}\Fdot$ in $\partial B_\epsilon$ and an open neighborhood $\Cal Y$ of $0$ in $\Bbb C$ such that,
for all $(\bold x, t)\in \Cal W^\prime\times \Cal Y$, for all $(a,b)\in\Bbb R^2-\{\bold 0\}$, $$\big(\bold x,\ a\,d_\bold
x(\operatorname{Re} (f+tL))+b\,d_\bold x(\operatorname{Im} (f+tL))\big)\ \not\in\ SS(\Fdot)+T_{{}_{\partial B_\epsilon}}^*X.\tag{$\ddag$}$$
Let $\Cal W:=\Cal W^\prime\cup(\partial B_\epsilon-\operatorname{supp}\Fdot)$. As $SS(\Fdot)+T_{{}_{\partial B_\epsilon}}^*X$ lies over
$\operatorname{supp}\Fdot\cap\partial B_\epsilon$, $(\ddag)$ also holds for all $(\bold x, t)\in \Cal W\times \Cal Y$ and
$(a,b)\in\Bbb R^2-\{\bold 0\}$.

\vskip .1in

 By $(\dagger)$ and Proposition 4.12, $SS(j_*j^*\Fdot)\subseteq SS(\Fdot)+N^*(B_\epsilon)$. In addition, if $\bold x\in\partial
B_\epsilon$ and $(\bold x, \omega)\in SS(\Fdot)+N^*(B_\epsilon)$, then certainly $(\bold x, \omega)\in SS(\Fdot)+T_{{}_{\partial B_\epsilon}}^*X$.
The desired conclusion now follows from $(\ddag)$.
\qed

\vskip .5in

Finally, we can prove

\vskip .3in

\noindent{\bf Theorem 5.3}. {\it  Suppose
that $\bold 0$ is an isolated point in the support of $\phi_f\Fdot$.

Then,  for every $\Fdot$-visible stratum $S_\alpha$, $(\bold 0, d_\bold 0 f)$ is an isolated point in
$\overline{T^*_{{}_{S_\alpha}}X}\cap\operatorname{im} d f$, and if we let 
$$
k_\alpha:=\left(\overline{T^*_{{}_{S_\alpha}}X}\ \cdot\ \operatorname{im} d f\right)_{(\bold 0, d_\bold 0 f)},
$$
then 
 $$ H^{i-1}(\phi_f\Fdot)_\bold 0\ \cong \ \bigoplus\Sb\Fdot{\text-visible}\\ S_\alpha\endSb\big(R^{{k_\alpha}}\otimes_{{}_R}\Bbb
H^{i-d_\alpha}(\Bbb N_\alpha,
\Bbb L_\alpha \ ; \ \Fdot)\big).$$ 

}

\vskip .3in

\noindent{\it Proof}. Fix a linear form $L$ and an $\epsilon>0$ so that Lemmas 5.1 and 5.2 hold; let $\Cal W$ and $\Cal Y$ be as in Lemma 5.2. 

If
necessary, select
$\epsilon$ smaller, so that all stratified critical points of $f$ inside $B_\epsilon$ occur on $V(f)$. Also choose
$\epsilon$ small enough, and a real $\omega_0>0$ small enough so that, if $0<\omega\leqslant\omega_0$,  $B_\epsilon\cap f^{-1}(\Bbb D_\omega)$ has
hypercohomogy isomorphic to the stalk cohomology of
$\Fdot$ at $\bold 0$ and so that, for all $v\in \Bbb D^*_{\omega_0}$,
$B_\epsilon\cap V(f-v)$ has the hypercohomology of the Milnor fibre. Also, choose $\omega_0$ small enough so that $\partial B_\epsilon\cap
f^{-1}(\Bbb D_{\omega_0})\subseteq \Cal W$. Fix a positive
$\omega<\omega_0$.

Let $\Bbb D_\omega(v)$ denote the closed disk
of radius
$\omega$, centered at $v$, and fix $v$ so small that $B_\epsilon\cap f^{-1}(\Bbb D_\omega)$ has hypercohomogy isomorphic to that of $B_\epsilon\cap
f^{-1}(\Bbb D_\omega(v))$; this is possible because all stratified critical points of $f$ occur on $V(f)$.  Again, by stratified Morse
theory, for all small
$t\in\Cal Y$, if
$\Psi_{L,t,v}(\bold z):=|f(\bold z)+tL(\bold z)-v|^2$, then $\Psi_{L,t,v}^{-1}[0, \omega^2]$ has the same hypercohomlogy as
$B_\epsilon\cap f^{-1}(\Bbb D_\omega(v))$, and $\Psi_{L,t,v}^{-1}(0)$ has the same hypercohomlogy as
$B_\epsilon\cap V(f-v)$. Now, apply microlocal Morse theory to $\Psi_{L,t,v}$ as its value goes from $0$ to $\omega^2$. 

By Lemma 4.10, there is no change in hypercohomology as $\Psi_{L,t,v}$ goes from $0$ to a sufficiently small positive value. After this,
Lemma 4.9, combined with Lemma 5.2, guarantees that stratified critical points on $\partial B_\epsilon$ produce no
change in hypercohomology. Hence, by Lemma 5.1, the remainder of the proof is precisely the same as that of Theorem 1.1.\qed

\enddocument